\documentclass[11pt,leqno]{article}
\usepackage[active]{srcltx}

\usepackage{amsfonts,latexsym,amsmath,amssymb,amsthm}
\usepackage{fullpage}
\newtheorem{theorem}{Theorem}[section]
\newtheorem{lemma}[theorem]{Lemma}

\newtheorem{proposition}[theorem]{Proposition}

\newtheorem{remark}[theorem]{Remark}

\theoremstyle{definition}

\theoremstyle{remark}

\newtheorem*{note*}{Note}

\numberwithin{equation}{section}

\makeatletter
\newcommand{\rank}{\mathop{\operator@font rank}}
\newcommand{\conv}{\mathop{\operator@font conv}}
\newcommand{\vol}{\mathop{\operator@font vol}}
\newcommand{\onetagright}{\tagsleft@false}
\makeatother

\newcommand{\ls}{\leqslant}
\newcommand{\gr}{\geqslant}
\renewcommand{\epsilon}{\varepsilon}
\newcommand{\prend}{$\quad \hfill \Box$}
\usepackage{color}

\usepackage{hyperref}

\begin{document}
\small

\title{\bf Random approximation and the vertex index of convex bodies}

\medskip

\author{Silouanos Brazitikos, Giorgos Chasapis and Labrini Hioni}

\date{}

\maketitle

\begin{abstract}
\footnotesize We prove that there exists an absolute constant $\alpha >1$ with the following
property: if $K$ is a convex body in ${\mathbb R}^n$ whose center of mass is at the origin,
then a random subset $X\subset K$ of cardinality ${\rm card}(X)=\lceil\alpha n\rceil $ satisfies with probability greater than $1-e^{-n}$
\begin{equation*}
K\subseteq c_1n\,{\rm conv}(X),
\end{equation*}
where $c_1>0$ is an absolute constant. As an application we show that the vertex index of any
convex body $K$ in ${\mathbb R}^n$ is bounded by $c_2n^2$, where $c_2>0$ is an absolute constant,
thus extending an estimate of Bezdek and Litvak for the symmetric case.
\end{abstract}

\section{Introduction}

The starting point of this article is the following result of Barvinok from \cite{Barvinok-2014}: If $C\subset\mathbb{R}^n$ is
a compact set then, for every $d>1$ there exists a subset $X\subseteq C$ of
cardinality ${\rm card}(X)\ls dn$ such that for any $z\in {\mathbb R}^n$ we have
\begin{equation}\label{eq:barvinok-1}
\max_{x\in X}|\langle z,x\rangle |\ls\max_{x\in C}|\langle z,x\rangle|\ls \gamma_d\sqrt{n}\max_{x\in X}|\langle z,x\rangle|,
\end{equation}
where $\gamma_d:=\frac{\sqrt{d}+1}{\sqrt{d}-1}$. For the proof of this fact, Barvinok assumes that the Euclidean
unit ball $B_2^n$ is the ellipsoid of minimal volume containing the convex hull ${\rm conv}(C)$ of $C$, and makes essential use of a
theorem of Batson, Spielman and Srivastava \cite{BSS-2009} on extracting an approximate John's decomposition with
few vectors from a John's decomposition of the identity. From \eqref{eq:barvinok-1} one can easily conclude that
if $K$ is an origin symmetric convex body in ${\mathbb R}^n$ then for any $d>1$ there exist $N\ls dn$ points $x_1,\ldots ,x_N\in K$
such that
\begin{equation}\label{eq:barvinok-2}{\rm absconv}(\{x_1,\ldots ,x_N\})\subseteq K\subseteq \gamma_d\sqrt{n}\,{\rm absconv}(\{x_1,\ldots ,x_N\}).\end{equation}
A generalization of Barvinok's lemma was recently obtained by the first named author in \cite{Brazitikos-diameter}:
There exists an absolute constant $\alpha >1$ with the following
property: if $K$ is a convex body whose minimal volume ellipsoid is the Euclidean unit ball,
then there exist $N\ls \alpha n$ points $x_1,\ldots ,x_N\in K\cap S^{n-1}$ such that
\begin{equation}\label{eq:brazitikos-1}
K\subseteq B_2^n\subseteq cn^{3/2}{\rm conv}(X),
\end{equation}
where $c>0$ is an absolute constant. The proof involves a more delicate theorem of Srivastava from \cite{Srivastava-2012}.
Using \eqref{eq:brazitikos-1} one can establish the following ``quantitative diameter version" of Helly's theorem (see \cite{Brazitikos-diameter}):
If $\{P_i: i\in I\}$ is a finite family of convex bodies in ${\mathbb R}^n$ with ${\rm diam}\left (\bigcap_{i\in I}P_i\right )=1$,
then there exist $s\ls \alpha n$ and $i_1,\ldots i_s\in I$ such that
\begin{equation}\label{eq:brazitikos-2}
{\rm diam}(P_{i_1}\cap\cdots\cap P_{i_s})\ls cn^{3/2},
\end{equation}
where $c>0$ is an absolute constant. Our first main result provides a random version of \eqref{eq:brazitikos-1} with an improved dependence
on the dimension.

\begin{theorem}\label{th:intro-1}There exists an absolute constant $\alpha >1$ with the following
property: if $K$ is a convex body in ${\mathbb R}^n$ whose center of mass is at the origin, if $N=\lceil\alpha n\rceil $ and if
$x_1,\ldots ,x_N$ are independent random points uniformly distributed in $K$ then, with probability greater than $1-e^{-n}$
we have
\begin{equation}\label{eq:intro-1}
K\subseteq c_1n\,{\rm conv}(\{ x_1,\ldots ,x_N\}),
\end{equation}
where $c_1>0$ is an absolute constant.
\end{theorem}

For the proof we may assume that $K$ is an isotropic convex body (see Section 2 for background information) and
we use the so-called one-sided $L_q$-centroid bodies of $K$; these are the convex bodies $Z_q^+(K)$, $q\gr 1$, with
support functions
\begin{equation}h_{Z_q^+(K)}(y)=\left ( 2\int_K\langle
x,y\rangle_+^qdx\right )^{1/q},\end{equation} where $a_+=\max
\{a,0\}$. We show that if $N\gr\alpha n$, where $\alpha >1$ is an absolute constant, then
$N$ independent random points $x_1,\ldots ,x_N$ uniformly distributed in $K$ satisfy
\begin{equation}{\rm conv}(\{x_1,\ldots ,x_N\})\supseteq c_1Z_2^+(K)\supseteq c_2 L_KB_2^n\end{equation} with
probability greater than $1-\exp (-n)$, where $c_1,c_2>0$ are absolute constants. Since $K$ is contained
in $(n+1)L_KB_2^n$, Theorem \ref{th:intro-1} follows.

We were led to Theorem \ref{th:intro-1} by the question to estimate the vertex index of a not necessarily symmetric $n$-dimensional
convex body. The vertex index of a symmetric convex body $K$ in ${\mathbb R}^n$ was introduced in \cite{Bezdek-Litvak-2007}
as follows:
\begin{equation}{\rm vi}(K)=\inf\Big\{ \sum_{j=1}^N\|y_j\|_K: K\subseteq {\rm conv}(\{ y_1,\ldots ,y_N\})\Big\},\end{equation}
where $\|\cdot\|_K$ is the norm with unit ball $K$ in ${\mathbb R}^n$. This index is closely related to the illumination
parameter of a convex body and to a well-known conjecture of Boltyanski and Hadwiger about covering of an $n$-dimensional convex body by
$2^n$ smaller positively homothetic copies (see \cite{Bezdek-Litvak-2007} and \cite{Gluskin-Litvak-2008}). Bezdek and Litvak proved that
\begin{equation}\frac{c_1n^{3/2}}{{\rm ovr}(K)}\ls {\rm vi}(K)\ls c_2n^{3/2},\end{equation}
where $c_1,c_2>0$ are absolute constants and ${\rm ovr}(K)$ is the outer volume ratio of $K$ (see Section 2 for the definition).
To the best of our knowledge the notion of vertex index has not been studied in the not necessarily symmetric case.
A way to define it for an arbitrary convex body $K$ in ${\mathbb R}^n$ is to consider first any $z\in {\rm int}(K)$ and
to set
\begin{equation}{\rm vi}_z(K)=\inf\Big\{\sum_{j=1}^Np_{K,z}(y_j): K\subseteq {\rm conv}(\{ y_1,\ldots ,y_N\})\Big\},\end{equation}
where
\begin{equation}p_{K,z}(x)=p_{K-z}(x)=\inf\{ t>0: x\in t(K-z)\}\end{equation}
is the Minkowski functional of $K$ with respect to $z$. Then, one may define the (generalized) vertex index of $K$ by
\begin{equation}{\rm vi}(K)={\rm vi}_{{\rm bar}(K)}(K),\end{equation}
where ${\rm bar}(K)$ is the center of mass of $K$. With this definition, we clearly have
${\rm vi}(K)={\rm vi}(K-{\rm bar}(K))$, and hence we may restrict our attention to centered convex bodies (i.e. convex
bodies whose center of mass is at the origin). In Section 4 we establish some elementary properties of this index and using Theorem \ref{th:intro-1}
we obtain the following general estimate.

\begin{theorem}\label{th:intro-3}There exist two absolute constants $c_1,c_2>0$ such that for every $n\gr 2$ and for every centered convex body $K$ in ${\mathbb R}^n$,
\begin{equation}\frac{c_1n^{3/2}}{{\rm ovr}({\rm conv}(K,-K))}\ls {\rm vi}(K)\ls c_2n^2.\end{equation}
\end{theorem}

A natural question, which is closely related to Theorem \ref{th:intro-1}, is to fix $N\gr\alpha n$ and to ask
for the largest value $t(N,n)$ for which $N$ independent random points $x_1,\ldots ,x_N$ uniformly distributed in $K$ satisfy
\begin{equation}{\rm conv}(\{x_1,\ldots ,x_N\})\supseteq t(N,n)\,K\end{equation} with
probability ``exponentially close" to $1$. A sharp answer to this question would unify Theorem \ref{th:intro-1} and the
following result from \cite{GM-concentration} which deals with the case where $N$ is exponential in $n$: For every $\delta \in (0,1)$
there exists $n_0=n_0(\delta )$ such that if $n\gr n_0$, if $C\log n/n\ls\gamma\ls 1$ and if $K$ is a centered
convex body in ${\mathbb R}^n$, then $N=\exp (\gamma n)$ independent random
points $x_1,\ldots ,x_N$ chosen uniformly from $K$ satisfy with
probability greater than $1-\delta $
\begin{equation}K\supseteq {\rm conv}(\{x_1,\ldots ,x_N\})\supseteq c(\delta )\gamma K,\end{equation}
where $c(\delta )$ is a constant depending on $\delta $. We prove the following.

\begin{theorem}\label{th:intro-2}Let $\beta\in (0,1)$. There exist a constant $\alpha =\alpha (\beta )>1$ depending only on $\beta $
and an absolute constant $c_1>0$ with the following
property: if $K$ is a centered convex body in ${\mathbb R}^n$, if $\alpha n\ls N\ls e^n$ and if
$x_1,\ldots ,x_N$ are independent random points uniformly distributed in $K$, then
\begin{equation}\label{eq:intro-2}
{\rm conv}(\{ x_1,\ldots ,x_N\})\supseteq \frac{c_1\beta\log (N/n)}{n}\,K.
\end{equation}
with probability greater than $1-e^{-N^{1-\beta }n^{\beta }}$.
\end{theorem}

In fact, Theorem \ref{th:intro-1} is a special case of Theorem \ref{th:intro-2}. The proof of both theorems is given
in Section~3.

\section{Notation and background}

We work in ${\mathbb R}^n$, which is equipped with a Euclidean structure $\langle\cdot ,\cdot\rangle $. We denote by $\|\cdot \|_2$
the corresponding Euclidean norm, and write $B_2^n$ for the Euclidean unit ball and $S^{n-1}$ for the unit sphere.
Volume is denoted by $|\cdot |$. We use the same notation $|X|$ for the cardinality of a finite set $X$.
We write $\omega_n$ for the volume of $B_2^n$ and $\sigma $ for the rotationally invariant probability
measure on $S^{n-1}$.

The letters $c,c^{\prime }, c_1, c_2,\ldots $ denote absolute positive constants which may change from line to line. Whenever we write
$a\simeq b$, we mean that there exist absolute constants $c_1,c_2>0$ such that $c_1a\ls b\ls c_2a$.  Also, if $K,L\subseteq \mathbb R^n$
we will write $K\simeq L$ if there exist absolute constants $c_1,c_2>0$ such that $ c_{1}K\subseteq L \subseteq c_{2}K$.

We refer to the book of Schneider \cite{Schneider-book} for basic facts from the Brunn-Minkowski theory and to the book
of Artstein-Avidan, Giannopoulos and V. Milman \cite{AGA-book} for basic facts from asymptotic convex geometry.

\smallskip

A convex body in ${\mathbb R}^n$ is a compact convex subset $K$ of ${\mathbb R}^n$ with non-empty interior. We say that $K$ is
symmetric if $x\in K$ implies that $-x\in K$, and that $K$ is centered if its center of mass
\begin{equation}{\rm bar}(K)=\frac{1}{|K|}\int_Kx\,dx \end{equation} is at the origin. The circumradius of $K$ is the radius of the
smallest ball which is centered at the origin and contains $K$:
\begin{equation}R(K)=\max\{ \|x\|_2:x\in K\}.\end{equation}
If $0\in {\rm int}(K)$ then the polar body $K^{\circ }$ of $K$ is defined by
\begin{equation}K^{\circ }:=\{ y\in {\mathbb R}^n: \langle x,y\rangle \ls 1
\;\hbox{for all}\; x\in K\}, \end{equation}
and the Minkowski functional of $K$ is defined by
\begin{equation}p_K(x)=\inf\{ t> 0: x\in tK\}.\end{equation}
Recall that $p_K$ is subadditive and positively homogeneous.

\smallskip

We say that a convex body $K$ is in John's position if the ellipsoid of maximal volume inscribed in $K$ is
the Euclidean unit ball $B_2^n$. John's theorem (see \cite[Chapter 2]{AGA-book}) states that
$K$ is in John's position if and only if $B_2^n\subseteq K$ and there exist $v_1,\ldots ,v_m\in {\rm
bd}(K)\cap S^{n-1}$ (contact points of $K$ and $B_2^n$) and positive real numbers $a_1,\ldots ,a_m$ such that
\begin{equation}\label{eq:bar-0}\sum_{j=1}^ma_jv_j=0\end{equation}
and the identity operator $I_n$ is decomposed in the form
\begin{equation}\label{eq:decomposition}I_n=\sum_{j=1}^ma_jv_j\otimes v_j, \end{equation}
where $(v_j\otimes v_j)(y)=\langle v_j,y\rangle v_j$. We say that a convex body $K$ is in L\"{o}wner's
position if the ellipsoid of minimal volume containing $K$ is the Euclidean unit ball
$B_2^n$. One can check that this holds true if and only if $K^{\circ }$ is in John's position;
in particular, we have a decomposition of the identity similar to
\eqref{eq:decomposition}. The outer volume ratio of a convex body $K$ in ${\mathbb R}^n$ is the quantity
\begin{equation}{\rm ovr}(K)=\inf\left\{ \left(\frac{|{\cal E}|}{|K|}\right )^{1/n}:{\cal E}\;\hbox{is an ellipsoid and}\;K\subseteq {\cal E}\right\}.\end{equation}
If $K$ is in L\"{o}wner's position then $(|B_2^n|/|K|)^{1/n}={\rm ovr}(K)$.

A convex body $K$ in ${\mathbb R}^n$ is called isotropic if it has volume $1$, it is centered,
and its inertia matrix is a multiple of the identity matrix:
there exists a constant $L_K >0$ such that
\begin{equation}\label{isotropic-condition}\int_K\langle x,\theta\rangle^2dx =L_K^2\end{equation}
for every $\theta $ in the Euclidean unit sphere $S^{n-1}$. It is known that if $K$ is isotropic then
\begin{equation}\label{eq:inclusions}cL_K\,B_2^n\subseteq K\subseteq (n+1)L_K\,B_2^n,\end{equation}
where $c>0$ is an absolute constant. The hyperplane conjecture
asks if there exists an absolute constant $C>0$ such that
\begin{equation}\label{HypCon}L_n:= \max\{ L_K:K\ \hbox{is isotropic in}\ {\mathbb R}^n\}\ls C\end{equation}
for all $n\gr 1$. Bourgain proved in \cite{Bourgain-1991} that $L_n\ls
c\sqrt[4]{n}\log\! n$, while Klartag \cite{Klartag-2006}
obtained the bound $L_n\ls c\sqrt[4]{n}$. A second proof of Klartag's bound
appears in \cite{Klartag-EMilman-2012}. We refer the reader to the article of V. Milman and Pajor \cite{VMilman-Pajor-1989} and to
the book \cite{BGVV-book} for an updated exposition of isotropic log-concave
measures and more information on the hyperplane conjecture.

The $L_q$-centroid body $Z_q(K)$ of $K$ is the
centrally symmetric convex body with support function
\begin{equation}\label{Zq-def}h_{Z_q(K)}(y)=\left(\int_K |\langle x,y\rangle|^{q}dx \right)^{1/q}.\end{equation}
Note that $K$ is isotropic if and only if it is centered and $Z_{2}(K)=
L_{K}B_2^n$. Also, if $T\in SL(n)$ then $Z_{q}(T(K))= T(Z_{q}(K))$. From H\"{o}lder's inequality it follows that
$Z_1(K)\subseteq Z_p(K)\subseteq Z_q(K)\subseteq Z_{\infty }(K)$ for
all $1\ls p\ls q\ls \infty $, where $Z_{\infty }(K)={\rm conv}(K,-K)$.
Using Borell's lemma (see \cite[Chapter 1]{BGVV-book}) one can check that
\begin{equation}\label{eq:Zq-inclusions} Z_q(K)\subseteq \overline{c}_1\frac{q}{p}Z_p(K)\end{equation}
for all $1\ls p<q$, where $\overline{c}_1>0$ is an absolute constant. In particular, if $K$ is isotropic then
\begin{equation}\label{eq:Zq-radius}R(Z_q(K))\ls \overline{c}_1qL_K.\end{equation} One can also check that if $K$ is
centered, then $Z_q(K)\supseteq c_2Z_{\infty }(K)$ for all $q\gr n$. For a proof of all
these assertions see \cite[Chapter 5]{BGVV-book}. The class of $L_q$-centroid bodies
of $K$ was introduced (with a different normalization) by Lutwak, Yang and Zhang in \cite{Lutwak-Yang-Zhang-2000}. An asymptotic
approach to this family was developed by Paouris in \cite{Paouris-GAFA} and \cite{Paouris-TAMS}.

For the proof of Theorem \ref{th:intro-2}
we generalize the arguments from \cite{Dafnis-Giannopoulos-Tsolomitis-2009} who used $L_q$-centroid bodies in order to
describe the asymptotic shape of the absolute convex hull of $N$ random points chosen from a convex body. The use
of one-sided $L_q$-centroid bodies allows one to consider the convex hull itself.

\section{Random approximation of convex bodies}

Let $K$ be a centered convex body of volume $1$ in ${\mathbb R}^n$.
For every $q\gr 1$ we consider the one-sided $L_q$-centroid body
$Z_q^+(K)$ of $K$ with support function
\begin{equation}h_{Z_q^+(K)}(y)=\left ( 2\int_K\langle
x,y\rangle_+^qdx\right )^{1/q},\end{equation} where $a_+=\max
\{a,0\}$. When $K$ is symmetric, it is clear that $Z_q^+(K)=Z_q(K)$. In any
case, we easily verify that
\begin{equation}Z_q^+(K)\subseteq 2^{1/q}Z_q(K).\end{equation}
Note that $Z_q^+(K)\subseteq 2^{1/q}K$ for all $q\gr 1$. Using Gr\"{u}nbaum's lemma (see \cite[Proposition 1.5.16]{AGA-book})
one can check that if $1\ls q\ls r<\infty $ then
\begin{equation}\label{eq:comp-q-r}\left (\frac{2}{e}\right )^{\frac{1}{q}-\frac{1}{r}}Z_q^+(K)\subseteq Z_r^+(K)\subseteq
\frac{Cr}{q}\left (\frac{2e-2}{e}\right )^{\frac{1}{q}-\frac{1}{r}}Z_q^+(K),\end{equation}
where $C>0$ is an absolute constant. The next lemma is due to Gu\'{e}don and E. Milman (see \cite{Guedon-EMilman-2011}).

\begin{lemma}\label{lem:two-norm}There exists an absolute constant $\overline{c}_0>0$ such that, for every isotropic convex body $K$ in ${\mathbb R}^n$,
\begin{equation}Z_2^+(K)\supseteq \overline{c}_0L_KB_2^n.\end{equation}
Equivalently, for any $\theta\in S^{n-1}$,
\begin{equation}h_{Z_2^+(K)}(\theta )=\left ( 2\int_K\langle
x,y\rangle_+^2dx\right )^{1/2}\gr \overline{c}_0L_K.\end{equation}
\end{lemma}

We also need the next lemma, which appears in \cite{Guedon-EMilman-2011} (see also \cite[Theorem 13.2.7]{BGVV-book}).

\begin{lemma}\label{lem:Zq+3}Let $K$ be a centered convex body of volume $1$ in ${\mathbb R}^n$.
We fix $\theta\in S^{n-1}$ and define $f_{\theta }(t)=|K\cap\{ x:\langle x,\theta\rangle =t\}|$, $t\in {\mathbb R}$. Then,
\begin{equation}\left (\frac{2}{e^2}\right )^{1/q}\left (\frac{\Gamma (n)\Gamma
(q+1)}{\Gamma (n+q+1)}\right )^{1/q}h_K(\theta )\ls
h_{Z_q^+(K)}(\theta )\ls
2^{1/q}h_K(\theta ).\end{equation}
\end{lemma}

\noindent {\it Proof.} We sketch the proof of the left hand side inequality. Let \begin{equation}H_{\theta }^+=\{ x\in {\mathbb R}^n:\langle
x,\theta\rangle\gr 0\}.\end{equation}First observe that, by the Brunn-Minkowski inequality, $f_{\theta }^{\frac{1}{n-1}}$ is concave
on its support, and hence we have
\begin{equation}f_{\theta }(t)\gr \Big ( 1-\frac{t}{h_K(\theta )}\Big )^{n-1}f_{\theta }(0)\end{equation}
for all $t\in [0,h_K(\theta )]$. Therefore,
\begin{align}
h^q_{Z_q^+(K)}(\theta ) = & 2\int_0^{h_K(\theta)}t^qf_{\theta }(t)dt
\gr 2\int_{0}^{h_K(\theta )}t^q\Big ( 1-\frac{t}{h_K(\theta )}\Big
)^{n-1}f_{\theta }(0)dt
\\ \nonumber
&= 2f_{\theta }(0)h_K^{q+1}(\theta )\int_{0}^1s^q(1-s)^{n-1}ds\\
\nonumber &= \frac{\Gamma (n)\Gamma (q+1)}{\Gamma (q+n+1)}\,2f_{\theta
}(0)h_K^{q+1}(\theta ).
\end{align}
Observe that
\begin{equation}2f_{\theta }(0)h_K(\theta )=\frac{f_{\theta }(0)}{\| f_{\theta }\|_{\infty }}\,2\|f_{\theta }\|_{\infty }h_K(\theta )\gr
\frac{f_{\theta }(0)}{\| f_{\theta }\|_{\infty }}\,(2|K\cap H_{\theta
}^+|).\end{equation}
We know that $\|f_{\theta }\|_{\infty }\ls e f_{\theta }(0)$ by a result of Fradelizi (see e.g. \cite[Theorem 2.2.2]{BGVV-book})
and that $|K\cap H_{\theta }^+|\gr e^{-1}$ by Gr\"{u}nbaum's lemma (see \cite[Proposition 1.5.16]{AGA-book}). Combining the above we get the result. \prend

\medskip

Theorem \ref{th:intro-2} and Theorem \ref{th:intro-1} will follow from the next fact, which generalizes work of Dafnis, Giannopoulos
and Tsolomitis \cite{Dafnis-Giannopoulos-Tsolomitis-2009} to the not necessarily symmetric setting.

\begin{theorem}\label{th:intro-4}Let $\beta\in (0,1)$. There exist a constant $\alpha =\alpha (\beta )>1$ depending only on $\beta $ 
and absolute constants $c_1,c_2>0$ with the following
property: if $K$ is a centered convex body in ${\mathbb R}^n$, if $N\gr\alpha n$ and if
$x_1,\ldots ,x_N$ are independent random points uniformly distributed in $K$ then there exists $q\gr c_1\beta\log (N/n)$ such that
\begin{equation}\label{eq:intro-4}
{\rm conv}(\{ x_1,\ldots ,x_N\})\supseteq c_2 Z_q^+(K)
\end{equation}
with probability greater than $1-e^{-N^{1-\beta }n^{\beta }}$.
\end{theorem}

Our proof of \eqref{eq:intro-4} is using the family of one-sided $L_q$-centroid bodies of $K$. In particular, we need
the following estimate.

\begin{lemma}\label{lem:correct-q}There exists an absolute constant $C>1$ with the following property: for every $n\gr 1$,
for every centered convex body $K$ in ${\mathbb R}^n$ and for every $q\gr 2$,
\begin{equation}\inf_{\theta\in S^{n-1}}\,\mu_K\left(\left\{ x:\langle x,\theta\rangle > \tfrac{1}{2} h_{Z_q^+(K)}(\theta )\right\}\right)\gr C^{-q} .\end{equation}
\end{lemma}

\noindent {\it Proof.} Let $K$ be a centered convex body in ${\mathbb R}^n$, let $q\gr 2$ and let $\theta\in S^{n-1}$.
We apply the Paley-Zygmund inequality \begin{equation}\label{eq:13.2.13}{\mathbb P}\,\big(g\gr t{\mathbb E}\,(g)\big)\gr
(1-t)^2\frac{[{\mathbb E}\,(g)]^2}{{\mathbb E}\,(g^2)}
\end{equation} for the non-negative random variable
\begin{equation}g_{\theta }(x)=2\langle x,\theta\rangle_+^q\end{equation} on $(K,\mu_K)$, where $\mu_K$ is Lebesgue
measure on $K$. Applying \eqref{eq:comp-q-r} with $r=2q$ we see that
\begin{equation}\label{eq:13.2.144}{\mathbb E}\,(g_{\theta }^2)=h^{2q}_{Z_{2q}^+(K)}(\theta )\ls C_1^q h^{2q}_{Z_q^+(K)}(\theta )=C_1^q\,[{\mathbb E}\,(g_{\theta })]^2,
\end{equation}where $C_1>0$ is an absolute constant. From \eqref{eq:13.2.13} we get
\begin{align} \mu_K(\{ x:\langle x,\theta\rangle > t\,h_{Z_q^+(K)}(\theta )\}) &= \mu_K(\{ x:\langle x,\theta\rangle > t\,[{\mathbb E}\,(g_{\theta })]^{1/q}\})
= \mu_K(\{ x:\langle x,\theta\rangle_+ > t\,[{\mathbb E}\,(g_{\theta })]^{1/q}\})\\
\nonumber &= \mu_K(\{ x:\langle x,\theta\rangle_+^q > t^q\,{\mathbb E}\,(g_{\theta })\}) =\mu_K(\{ x:g_{\theta }(x) > 2t^q\,{\mathbb E}\,(g_{\theta })\})\\
\nonumber &\gr (1-2t^q)^2\frac{[{\mathbb E}\,(g_{\theta })]^2}{{\mathbb E}\,(g_{\theta }^2)}\gr \frac{(1-2t^q)^2}{C_1^q}
\end{align}
for every $t\in (0,2^{-\frac{1}{q}})$. Choosing $t=\tfrac{1}{2}$ we get the lemma with $C =4C_1$.\prend

\bigskip

\noindent {\bf Proof of Theorem \ref{th:intro-4}}. 
Let $q\gr 2$ and consider the random polytope $C_N:={\rm conv}\{x_1,\ldots ,x_N\}$.
With probability equal to one, $C_N$ has non-empty interior and, for
every $J=\{ j_1,\ldots ,j_n\}\subset \{ 1,\ldots ,N\}$, the points
$x_{j_1},\ldots ,x_{j_n}$ are affinely independent. Write $H_J$ for
the affine subspace determined by $x_{j_1},\ldots ,x_{j_n}$ and
$H_J^+$, $H_J^-$ for the two closed halfspaces whose bounding
hyperplane is $H_J$.

If $\tfrac{1}{2}Z_q^+(K)\not\subseteq C_N$, then there exists $x\in \tfrac{1}{2}Z_q^+(K)\setminus C_N$, and hence,
there is a facet of $C_N$ defining some affine subspace $H_J$ as above that satisfies the following: either $x\in
H_J^-$ and $C_N\subset H_J^+$, or $x\in H_J^+$ and $C_N\subset
H_J^-$. Observe that, for every $J$, the probability of each of these two events
is bounded by \begin{equation}\Big (\sup_{\theta\in S^{n-1}}\mu_K\Big(\Big\{ x:\langle x,\theta \rangle\ls \tfrac{1}{2}h_{Z_q^+(K)}(\theta )\Big\}\Big )\Big
)^{N-n}\ls\bigl (1-C^{-q}\bigr )^{N-n},\end{equation}where $C>0$ is the constant in Lemma \ref{lem:correct-q}. It follows that
\begin{equation}\label{6.24}{\mathbb P}\left (\tfrac{1}{2}Z_q^+(K)\not\subseteq C_N\right )\ls
2\binom{N}{n} (1-C^{-q})^{N-n}.\end{equation}
Since $\binom{N}{n}\ls\left (\frac{eN}{n}\right)^n$, this probability is smaller than $\exp (-N^{1-\beta }n^{\beta })$ if
\begin{equation}\left(\frac{2eN}{n}\right )^n(1-C^{-q})^{N-n}<\left(\frac{2eN}{n}\right )^ne^{-C^{-q}(N-n)}<\exp (-N^{1-\beta }n^{\beta }),\end{equation}
and the second inequality is satisfied if
\begin{equation}\label{eq:cond-q}\frac{N}{n}-1> C^q\left [\left (\frac{N}{n}\right )^{1-\beta }+\log\left(\frac{2eN}{n}\right )\right ].\end{equation}
We choose $q =\frac{\beta }{2\log C}\log\left (\frac{N}{n}\right )$ and $\alpha_1(\beta ):=C^{4/\beta }$. Note that if $N\gr\alpha_1(\beta ) n$ then $q\gr 2$ if and
that \eqref{eq:cond-q} becomes
\begin{equation}\label{eq:cond-q-2}\frac{N}{n}-1>\left (\frac{N}{n}\right )^{1-\frac{\beta }{2}}+\left (\frac{N}{n}\right )^{\frac{\beta }{2}}\log\left(\frac{2eN}{n}\right ).\end{equation}
Since 
\begin{equation}\lim_{t\to +\infty }\left [t-1-t^{1-\frac{\beta }{2}}-t^{\frac{\beta }{2}}\log (2et)\right ]=+\infty ,\end{equation}
we may find $\alpha_2(\beta )$ such that \eqref{eq:cond-q-2} is satisfied for all $N\gr \alpha_2(\beta ) n$.
Setting $\alpha =\max\{\alpha_1(\beta ),\alpha_2(\beta )\}$ we see that the assertion of the theorem is satisfied
with probability greater that $1-e^{-N^{1-\beta }n^{\beta }}$ for all $N\gr \alpha n$, with $q \gr c_2\beta \log\left (\frac{N}{n}\right )$,
where $c_2>0$ is an absolute constant. \prend

\bigskip

\noindent {\bf Proof of Theorem \ref{th:intro-2}.} Let $\beta\in (0,1)$ and let $\alpha =\alpha (\beta )$
be the constant from Theorem \ref{th:intro-4}. Let $\alpha n\ls N\ls e^n$ and let $x_1,\ldots ,x_N$ be independent random points uniformly
distributed in $K$. Applying Lemma \ref{lem:Zq+3} with $q=n$ we see that $h_{Z_n^+(K)} \gr c_1h_K(\theta )$ for all
$\theta\in S^{n-1}$, and hence
\begin{equation}Z_n^+(K)\supseteq c_1K,\end{equation}
where $c_1>0$ is an absolute constant. From Theorem \ref{th:intro-4} we know that if $q= c_2\beta\log (N/n)$ (note also that $q\ls n$) then 
\begin{equation}C_N={\rm conv}(\{x_1,\ldots ,x_N\})\supseteq c_3Z_q^+(K)\end{equation}
with probability greater than $1-\exp (-N^{1-\beta }n^{\beta })$, where $c_2,c_3>0$ are absolute constants. From \eqref{eq:comp-q-r}
we see that \begin{equation}Z_n^+(K)\subseteq
\frac{c_4n}{q}\left (\frac{2e-2}{e}\right )^{\frac{1}{q}-\frac{1}{n}}Z_q^+(K)\subseteq \frac{2c_4n}{q}Z_q^+(K),\end{equation}
where $c_4>0$ is an absolute constant. Combining the above we get that
\begin{equation}C_N={\rm conv}(\{x_1,\ldots ,x_N\})\supseteq \frac{c_5q}{n}\,K\supseteq \frac{c_6\beta\log (N/n)}{n}\,K\end{equation}
with probability greater than $1-\exp (-N^{1-\beta }n^{\beta })$, where $c_5,c_6>0$ are absolute constants. \prend

\medskip

Choosing $N=\lceil\alpha n\rceil $ in Theorem \ref{th:intro-2} we immediately get Theorem \ref{th:intro-1}:

\begin{theorem}\label{th:main-gen}There exists an absolute constant $\alpha >1$ with the following
property: if $K$ is a centered convex body in ${\mathbb R}^n$ then a random subset $X\subset K$ of
cardinality ${\rm card}(X)=\lceil\alpha n\rceil $ satisfies
\begin{equation}\label{eq:iso-1}
K\subseteq Cn\,{\rm conv}(X)
\end{equation}
with probability greater than $1-e^{-n}$, where $C>0$ is an absolute constant.
\end{theorem}

\section{Generalized vertex index}\label{sec:gvi}

Let $K$ be a convex body in ${\mathbb R}^n$. From the definition of the vertex index that we gave in the introduction,
we may clearly assume that $K$ is centered, and then
\begin{equation}{\rm vi}(K)=\inf\Big\{\sum_{j=1}^Np_{K}(y_j): K\subseteq {\rm conv}(\{ y_1,\ldots ,y_N\})\Big\},\end{equation}
where $p_K$ is the Minkowski functional of $K$. Since every origin symmetric convex body is centered, our definition coincides with
the one given by Bezdek and Litvak in \cite{Bezdek-Litvak-2007} for the symmetric case.

It is also easy to check that the vertex index is invariant under invertible linear transformations.
For every convex body $K$ in ${\mathbb R}^n$ and any $T\in GL(n)$ one has
\begin{equation}{\rm vi}(T(K))={\rm vi}(K).\end{equation}
To see this, note that $T(K)\subseteq {\rm conv}(\{ y_1,\ldots ,y_N\})$ if and only if $K\subseteq {\rm conv}(\{ x_1,\ldots ,x_N\})$ where
$T(x_j)=y_j$, therefore
\begin{align}{\rm vi}(T(K)) &=\inf\Big\{\sum_{j=1}^Np_{T(K)}(T(x_j)): K\subseteq {\rm conv}(\{ x_1,\ldots ,x_N\})\Big\}\\
\nonumber &=\inf\Big\{\sum_{j=1}^Np_{K}(x_j): K\subseteq {\rm conv}(\{ x_1,\ldots ,x_N\})\Big\}\\
\nonumber &={\rm vi}(K).
\end{align}
Another useful observation is that the vertex index is stable under a variant of the Banach-Mazur distance. Recall that
the Banach-Mazur distance between two convex bodies $K$ and $L$ in ${\mathbb R}^n$ is the quantity
\begin{equation}d(K,L)=\inf\{ t>0: T(L+y)\subseteq K+x\subseteq
t(T(L+y))\},\end{equation} where the infimum is over all $T\in GL(n)$ and $x,y\in {\mathbb R}^n$.
Given two centered convex bodies $K$ and $L$, we set
\begin{equation}\tilde{d}(K,L)=\inf\{ t>0: T(L)\subseteq K\subseteq
tT(L)\},\end{equation} where the infimum is over all $T\in GL(n)$. Note that if $K$ and $L$ are
symmetric convex bodies then $\tilde{d}(K,L)=d(K,L)$. With this definition we easily check that
if $K$ and $L$ are centered convex bodies in ${\mathbb R}^n$ then
\begin{equation}{\rm vi}(K)\ls \tilde{d}(K,L)\,{\rm vi}(L).\end{equation}
The main result of this section is the upper bound in Theorem \ref{th:intro-3}.

\begin{proposition}\label{prop:intro-3-upper}There exists an absolute constant $C_1>0$ such that,
for every $n\gr 2$ and for every centered convex body $K$ in ${\mathbb R}^n$,
\begin{equation}{\rm vi}(K)\ls C_1n^2.\end{equation}
\end{proposition}

\noindent {\it Proof.} We may assume that $K$ is isotropic. By Theorem \ref{th:main-gen} we can find $N\ls\alpha n$ and $x_1,\ldots ,x_N\in K$
such that
\begin{equation}
K\subseteq Cn\,{\rm conv}(\{x_1,\ldots ,x_N\}),
\end{equation}
where $\alpha ,C>0$ are absolute constants. We set $y_j=Cnx_j$, $1\ls j\ls N$. Then, $K\subseteq {\rm conv}(\{y_1,\ldots ,y_N\})$ and $p_K(y_j)=Cnp_K(x_j)\ls Cn$, therefore
\begin{equation}{\rm vi}(K)\ls\sum_{j=1}^Np_K(y_j)\ls CnN\ls C\alpha n^2.\end{equation}
The result follows with $C_1=C\alpha $. \prend

\medskip

For the lower bound we just check that the argument of \cite{Bezdek-Litvak-2007} remains valid in the not necessarily
symmetric case.

\begin{proposition}\label{prop:intro-3-lower}There exists an absolute constant $c>0$ such that,
for every $n\gr 2$ and for every centered convex body $K$ in ${\mathbb R}^n$,
\begin{equation}{\rm vi}(K)\gr\frac{cn^{3/2}}{{\rm ovr}({\rm conv}(K,-K))}.\end{equation}
\end{proposition}

\noindent {\it Proof.} By the linear invariance of the vertex index we may assume that $B_2^n$ is the ellipsoid of minimal volume which contains ${\rm conv}(K,-K)$.
In other words, $K\subseteq {\rm conv}(K,-K)\subseteq B_2^n$ and
\begin{equation}\left (\frac{|B_2^n|}{|{\rm conv}(K,-K)|}\right )^{1/n}= {\rm ovr}({\rm conv}(K,-K)).\end{equation}
For any $N\in {\mathbb N}$ and $y_1,\ldots ,y_N$ such that $K\subseteq {\rm conv}(\{y_1,\ldots ,y_N\})$
we consider the absolute convex hull $Q={\rm conv}(\{\pm y_1,\ldots
,\pm y_N\})\supseteq {\rm conv}(K,-K)$ of $y_1,\ldots ,y_N$. Then,
\begin{equation}Q^{\circ }=\{ x\in {\mathbb R}^n:\,|\langle x,y_j\rangle |\ls
1\;\hbox{for all}\;j=1,\ldots ,N\},\end{equation} and a result of
Ball and Pajor \cite{Ball-Pajor-1990} provides the lower bound
\begin{equation}|Q^{\circ }|\gr\left (\frac{n}{\sum_{j=1}^N\| y_j\|_2}\right
)^{1/n}\end{equation} for its volume. Using the Blaschke-Santal\'{o} inequality we get
\begin{equation}|{\rm conv}(K,-K)|\ls |Q|\ls \frac{|B_2^n|^2}{|Q^{\circ }|}\ls |B_2^n|^2\left
(\frac{\sum_{j=1}^N\| y_j\|_2}{n}\right )^n.\end{equation} It
follows that \begin{equation}1 \ls  \left (\frac{|B_2^n|}{|{\rm conv}(K,-K)|}\right
)^{1/n}|B_2^n|^{1/n}\frac{\sum_{j=1}^N\| y_j\|_2}{n}\ls \frac{{\rm ovr}({\rm conv}(K,-K)))}{cn^{3/2}}\sum_{j=1}^N\| y_j\|_2\end{equation}
for some absolute constant $c>0$. Since $K\subseteq B_2^n$, we have $\| y_j\|_2\ls p_{K}(y_j)$ for all
$j=1,\ldots ,N$. Therefore,
\begin{equation}\sum_{j=1}^Np_{K}(y_j)\gr\frac{cn^{3/2}}{{\rm ovr}({\rm conv}(K,-K))},\end{equation}
and taking the infimum over all $N$ and all such $N$-tuples $(y_1,\ldots ,y_N)$ we get the lower bound for $
{\rm vi}(K)$.
\prend

\begin{remark}\label{rem:final}\rm The lower bound of Proposition \ref{prop:intro-3-lower} is not sharp, even in the symmetric case. Gluskin
and Litvak \cite{Gluskin-Litvak-2012} have proved that for every $n\gr 1$ there exists a symmetric convex body $K$ in ${\mathbb R}^n$
such that
\begin{equation}{\rm ovr}(K)\gr c\sqrt{\frac{n}{\log (2n)}}\quad\hbox{and}\quad {\rm vi}(K)\gr cn^{3/2}.\end{equation}
It would be interesting to provide alternative lower bounds for ${\rm vi}(K)$; in particular, it would be interesting to
decide whether, in the non-symmetric case, the upper bound ${\rm vi}(K)\ls Cn^2$ of Proposition \ref{prop:intro-3-upper} is sharp or not.
\end{remark}

\bigskip

\noindent {\bf Acknowledgement.} We would like to thank Apostolos Giannopoulos for useful discussions. The first named author acknowledges support from Onassis Foundation.

\bigskip

\bigskip

\footnotesize
\bibliographystyle{amsplain}

\bigskip

\bigskip

\thanks{\noindent {\bf Keywords:}  Convex bodies, isotropic position, centroid bodies, random polytopal approximation.}

\smallskip

\thanks{\noindent {\bf 2010 MSC:} Primary 52A23; Secondary 52A35, 46B06, 60D05.}

\bigskip

\bigskip

\noindent \textsc{Silouanos \ Brazitikos}: Department of
Mathematics, University of Athens, Panepistimioupolis 157-84,
Athens, Greece.

\smallskip

\noindent \textit{E-mail:} \texttt{silouanb@math.uoa.gr}

\bigskip

\noindent \textsc{Giorgos \ Chasapis}: Department of
Mathematics, University of Athens, Panepistimioupolis 157-84,
Athens, Greece.

\smallskip

\noindent \textit{E-mail:} \texttt{gchasapis@math.uoa.gr}

\bigskip

\noindent \textsc{Labrini \ Hioni}: Department of
Mathematics, University of Athens, Panepistimioupolis 157-84,
Athens, Greece.

\smallskip

\noindent \textit{E-mail:} \texttt{lamchioni@math.uoa.gr}

\bigskip

\end{document}